\documentclass[12pt]{article}

\usepackage{amssymb,amsmath,exscale}

\usepackage{graphicx}

\usepackage[applemac]{inputenc}
\usepackage[T1]{fontenc}
\usepackage{dsfont}

\usepackage{color}
\definecolor{marin}{rgb}   {0.,   0.3,   0.7} 
\definecolor{rouge}{rgb}   {0.8,   0.,   0.} 
\definecolor{sepia}{rgb}   {0.8,   0.5,   0.} 
\usepackage[colorlinks,citecolor=marin,linkcolor=rouge,
            bookmarksopen,
            bookmarksnumbered
           ]{hyperref}

\newtheorem{lemma}{Lemma}[section]

\newtheorem{theorem}[lemma]{Theorem}
\newtheorem{proposition}[lemma]{Proposition}

\newtheorem{remark}[lemma]{Remark}
\newtheorem{example}[lemma]{Example}

\newtheorem{notation}[lemma]{Notation}
\newtheorem{definition}[lemma]{Definition}
\newtheorem{conclusion}[lemma]{Conclusion}
\numberwithin{equation}{section}
\newcommand{\QED}{\mbox{}\hfill \raisebox{-0.2pt}{\rule{5.6pt}{6pt}\rule{0pt}{0pt}} 
          \medskip\par}             
\newenvironment{Proof}{\noindent
    \parindent=0pt\abovedisplayskip = 0.5\abovedisplayskip
    \belowdisplayskip=\abovedisplayskip{\bfseries Proof. }}{\QED}

\newcommand{\dd}{\mathrm{d}}
\newcommand{\E}{\mathbb{E}}

\newcommand{\jb}{{\boldsymbol{j}}}
\newcommand{\kb}{{\boldsymbol{k}}}

\newcommand{\N}{\mathbb{N}}

\newcommand{\R}{\mathbb{R}}

\newcommand{\Ccal}{\mathcal{C}}

\newcommand{\T}{\mathbb{T}}

\newcommand{\F}{\mathcal{F}}
\newcommand{\PR}{\mathbb{P}}

\newcommand{\ds}{\displaystyle}

\newcommand{\Norm}[2]{\|#1\|\left.\vphantom{T_{j_0}^0}\!\!\right._{#2}}         
\newcommand{\SNorm}[2]{|#1|\left.\vphantom{T_{j_0}^0}\!\!\right._{#2}}             

\title{Weak backward error analysis for SDEs}        
\author{Arnaud Debussche and Erwan Faou\\[4ex]
\small INRIA \& ENS Cachan Bretagne,  IRMAR-UMR 6625\\[-1ex]
\small Avenue Robert Schumann
F-35170 Bruz\\[-1ex]
\small\it email: \tt Arnaud.Debussche@bretagne.ens-cachan.fr, \\[-1ex]\small \tt Erwan.Faou@inria.fr\\[4ex]
}       

\begin{document}

\maketitle

\abstract{
We consider numerical approximations of stochastic differential equations by the Euler method. In the case where the SDE is elliptic or hypoelliptic, we show a weak backward error analysis result in the sense that the generator associated with the numerical solution coincides with the solution of a modified Kolmogorov equation up to high order terms with respect to the stepsize. This implies that every invariant measure of the numerical scheme is close to a modified invariant measure obtained by asymptotic expansion. Moreover, we prove that, up to negligible terms,  the dynamic 
associated with the Euler scheme is exponentially mixing. }
\bigskip

\begin{centerline}
{\bf R\'esum\'e}
\smallskip

Nous \'etudions la discr\'etisation d'une \'equation diff\'erentielle stochastique (EDS) par le sch\'ema d'Euler. Dans le cas d'une EDS elliptique ou hypoelliptique nous montrons un r\'esultat d'analyse 
d'erreur r\'etrograde : une fonctionnelle de la solution num\'erique est 
proche de  la solution d'une \'equation de Kolmogorov modifi\'ee \`a des ordres arbitrairement
\'elev\'es par rapport au pas de discr\`etisation.  On obtient ainsi que toute mesure invariante du sch\'ema num\'erique
est proche d'une mesure invariante modifi\'ee obtenue par d\'eveloppement asymptotique. 
De plus, le sch\'ema est exponentiellement m\'elangeant \`a des ordres arbitrairement
\'elev\'es.

\bigskip
\end{centerline}




\section{Introduction}

In the last decades, {\em backward error analysis} has become one of the most powerful tool to analyze the long time behavior of numerical schemes applied to evolution equations. The main idea can be described as follows: Let us consider an ordinary differential equation of the form 
$$
\dot y(t) = f(y(t)), 
$$
where $f: \R^n \to \R^n$ is a smooth vector field, and denote by $\varphi^f_t(y)$ the associated flow. By definition, a numerical method defines for a small time step $\tau$ an approximation $\Phi_\tau$ of the exact flow $\varphi_\tau$: We have for bounded $y \in \R^n$, $\Phi_\tau(y) = \varphi^f_\tau(y) + \mathcal{O}(\tau^{r+1})$ where $r$ is the order of the method. 

The idea of backward error analysis is to show that $\Phi_\tau$ can be interpreted as the exact flow $\varphi_\tau^{f_\tau}$ or a {\em modified vector field} defined as a series in powers of $\tau$
$$
f_\tau = f + \tau^{r} f_r + \tau^{r+1} f_{r+1} + \cdots , 
$$
where $f_\ell$, $\ell \geq r$ are vector fields depending on the numerical method. 
In general, the series defining $f_\tau$ does not converge, but it can be shown that for bounded $y$, we have for arbitrary $N$ 
$$
\Phi_\tau(y) = \varphi_\tau^{f_\tau^N}(y) + C_N \tau^{N},
$$
where $f_\tau^N$ is the truncated series:
$$
f_\tau^N = f + \tau^{r} f_r +  \cdots + \tau^N f_N. 
$$
Under some analyticity assumptions, the constant $C_N \tau^N$ can be optimized in $N$, so that the error term in the previous equation can be made exponentially small with respect to $\tau$. 

Such a result is very important and has many applications in the case where $f$ has some strong geometric properties, such as Hamiltonian or reversible structure. In this situation, and under some compatibility conditions on the numerical method $\Phi_\tau$,  the modified vector field $f_\tau$ inherits the structure of $f$. For example if $\Phi_\tau$ is symplectic and $f$ Hamiltonian, then $f_\tau$ remains Hamiltonian. This has major consequences such as the preservation of a modified Hamiltonian over very long time (of order $\tau^{-N}$) for the numerical solution, from which we can deduce long time stability results, existence of numerical invariant tori in the integrable case, etc... 

In the Hamiltonian case, this idea goes back to Moser \cite{Moser68}, but was applied later to symplectic integrator by Benettin \& Giorgilli \cite{BeGi94}, Hairer \& Lubich \cite{HL97} and Reich \cite{Reich99}. Such results now form the core of the modern {\em geometric numerical integration theory} for which we refer to the classical textbooks \cite{HLW} and \cite{Reic04}. 

More recently, these ideas have been extended in some situations to Hamiltonian PDEs: First in the linear case \cite{DF09}, and then in the semilinear case (nonlinear Schr\"odinger or wave equations), see \cite{FG,Fbook}. 

As far as stochastic differential equations (SDEs)
are concerned, this approach has not been developed very much so far. 
Let us recall that given a SDE in $\R^d$ of the form
$$
\dd X=f(X)\dd t+\sigma(X)\dd W
$$
discretized by a numerical scheme - such as the Euler scheme for instance - with time step 
$\tau$ providing a discrete sequence $(X_p)_{p\in \N}$, then the error can be measured in the strong 
or weak sense. Strong error means that $X_p$ is a pathwise approximation of 
$X(p\tau)$, and it is well known that the Euler scheme has strong order $1/2$. Under standard 
assumptions on $f$ and $\sigma$, we have
$$
\E\left(\sup_{p=0,\dots,[T/\tau]}\| X_p-X(t_p)\|^k\right)\le c_1(k,T) \tau^{1/2}, \quad k\ge 1,\; T>0.
$$
In this work, we consider another error which is often more important. We investigate the weak error which concerns the law of the solution. The Euler scheme has weak order $1$. 
Under suitable smoothness assumptions on 
$f$, $\sigma$ and $\varphi\;:\; \R^d\to \R$ (see for instance \cite{kloeden-platten, milstein-tretyakov, Talay96}):
$$
\left| \E\left(\varphi(X_p)\right)- \E\left(\varphi(X(t_p))\right)\right|\le c_2(\varphi, T)\tau,
\quad p=0,\dots, [T/\tau],\quad T>0.
$$ 
An attempt  
has been made by Shardlow \cite{Shardlow} to extend the backward error analysis to this 
context. He has  shown that the construction of a 
modified SDE associated with the Euler scheme can be performed, but only at the first step, {\it ie} for $N=2$,
and only for additive noise, {\it ie} when $\sigma(X)$ does not depend on $X$. In this case, he is able to write down a 
modified SDE:
$$
\dd\tilde X=\tilde f(\tilde X)\dd t+\tilde\sigma(\tilde X)\dd W
$$
such that 
$$
\left| \E\left(\varphi(X_p)\right)- \E\left(\varphi(\tilde X(t_p))\right)\right|\le  c_3(\varphi, T)\tau^2,
\quad p=0,\dots, [T/\tau],\quad T>0.
$$ 
He explains that for multiplicative noise or higher order, 
there are too many conditions to be satisfied by the coefficients of the modified equations. 

In this paper we take another approach, and build a modified equation not at the level of the SDE, but at the level of the generator associated with the process solution of the SDE. It is 
well known that given  $\varphi : \R^d\to \R$ and denoting by 
$X(t,x)$ the solution of the SDE satisfying $X(0)=x$, the function 
$u(t,x)=\E(\varphi(X(t,x)))$ satisfies the Kolmogorov equation 
$$
\partial_t u(t,x) = L(x,\partial_x) u(t,x), 
$$
where $L$ is the order $2$ Kolmogorov  operator associated with the SDE.
 
 In the case of the Euler method applied to a SDE, we show that with the numerical solution, we can associate {\em modified Kolmogorov operator} of the form 
$$
L(\tau,x,\partial_x) = L(x,\partial_x) + \tau L_1(x,\partial_x) + \tau L_2(x,\partial_x) + \cdots
$$
where $L_\ell$, $\ell\geq 1$ are some modified operator of order $2\ell + 2$. Again, the series
does not converge but truncated series:
$$
L^{(N)}(\tau,x,\partial_x) = L(x,\partial_x) + \tau L_1(x,\partial_x) + \cdots +\tau L_N(x,\partial_x) 
$$
are considered.

Note that in contrast with the classical case, we do not have a modified SDE and cannot straightforwardly define a solution to the {\em modified equation}
$$
\partial_t v^N(t,x) = L^{(N)}( \tau, x,\partial_x) v^N(t,x).
$$
However, in the case where the SDE is elliptic or hypoelliptic, we can build an approximated solution $v^{(N)}$  such that 
$$
\left| \E\left(\varphi(X_p)\right)- v^{(N)}(p\tau,x)\right|\le  c_4(\varphi, T,N)\tau^{N},
\quad p=0,\dots, [T/\tau],\quad T>0.
$$
Furthermore, using the exponential convergence to equilibrium, we prove that in fact the constant
$c_4$ does not depend on $T$ so that we have an approximation result valid on very long times.
We also show that there exist a modified invariant measure for $L^{(N)}(\tau,x,\partial_x)$. 

We can then use this weak backward error analysis to prove that the numerical solution $X_p$, $p \geq 1$ obtained by the Euler  scheme is exponentially mixing up to some very small error, and for all times. 
This is typically a geometric numerical integration result in the sense that we prove the persistence of a qualitative property of the exact  flow (exponential mixing) to the numerical approximation, over long times. 

Note that $v^{(N)}$ is in fact constructed as a truncated series:
$
v^{(N)}=\sum_{n=0}^N \tau^n v_n
$
and that $v_0=u$ is the solution of the Kolmogorov equation. Therefore, our result provides
an expansion of the error as in \cite{TalayTubaro} (see also 
\cite{BallyTalay1, BallyTalay2}).  However, the expansion is different here. 

Error estimates on long times for elliptic and hypoelliptic SDEs have already been proved. In \cite{Mattingly1, Talay90,Talay96,TalayLangevin},
it is shown that for a sufficiently small time step the Euler scheme defines an ergodic  process and that the invariant measure of the Euler scheme is close to the invariant
measure of the SDE. In \cite{TalayTubaro}, the first term of an expansion of the invariant measure
of the Euler scheme with respect to $\tau$ is also given. In our
work, we provide the expansion at any order.

We emphasize that in our result, there is no particular smallness assumption on the stepsize $\tau$ used to define the numerical solution. In particular, the discrete process is not supposed to have a unique invariant measure, as in \cite{Mattingly1} or \cite{Talay90,Talay96,TalayLangevin}. 

This is also the case in the recent work \cite{Mattingly2}. There, it is shown that given an 
elliptic or hypoelliptic SDE, the ergodic 
averages provided by the Euler scheme 
are asymptotically close to the average of the invariant measure of the
SDE.
Higher order schemes  also considered. The main tool in
\cite{Mattingly2} is the ellipticity or hypoellipticity of the Poisson equation, {\it ie} the  equation $L(x,\partial_x)u=g$.

As in \cite{Mattingly2}, we consider the case where the SDE is set on the torus $\T^d$. 
This simplifies the presentation and the main ideas are not hidden by technical 
difficulties. In the same spirit, we only study the Euler scheme. In a forthcoming article,
we will present more realistic applications of our method for  SDEs set on $\R^n$ with polynomial growth coefficients under suitable assumptions, and for more general schemes. As an example,
we will treat the Langevin equation as in \cite{TalayLangevin}. 

\section{Preliminaries}

We consider the stochastic differential equation 
$$
\dd X(t) = f(X(t)) \dd t + \sigma(X(t)) \dd W, 
$$
where the unknown $X=(X^i)_{i=1\dots,d}$ lives in the $d$-dimensional torus $\T^d$.  Also,
$f=(f^i)_{i = 1}^d$ and  $\sigma= (\sigma^i_\ell(x))_{i = 1,\dots,d, \ell=1,\dots,m}$ are smooth vector fields periodic in $x \in \T^d$. The process $(W^1(t),\ldots,W^m(t))$ is a $m$-dimensional standard Wiener process over a probability space $(\Omega,\F,\PR)$ endowed with a filtration $(\F_t)_{t\ge 0}$.
Using these notations, we rewrite the equation as:
$$
\dd X^i(t) = f^i(X(t)) \dd t + \sum_{\ell = 1}^m \sigma^i_\ell(X(t)) \dd W^\ell(t),\quad i=1,\dots, d. 
$$

In all the paper, smooth functions means $\Ccal^\infty $ functions. 
Given a smooth function $\psi$ defined on $\T^d$, we denote by $\Norm{\psi}{\Ccal^k}$
its norm in $\Ccal^k(\T^d,\R)$. We also denote by $\|\psi\|_\infty=\sup_{x\in\T^d}
|\psi(x)|$. 
For a multiindex $\kb = (k_1,\ldots,k_d) \in \N^d$, we set $|\kb| = k_1 + \cdots + k_d$ 
and
$$
\partial_\kb \psi(x) = \frac{\partial^{|\kb|}\psi(x)}{\partial x_1^{k_1} \cdots \partial x_d^{k_d}}, \quad x\in \T^d.
$$
Therefore 
$$
\Norm{\psi}{\Ccal^k} := \sup_{\substack{\jb = (j_1,\ldots,j_d) \\ |\jb| \leq k} } |  \partial_\kb \psi(x) |. 
$$
We also define the semi-norm
$$
\SNorm{\psi}{\Ccal^k} := \sup_{\substack{\jb = (j_1,\ldots,j_d) \\ 1\le |\jb| \leq k} } |  \partial_\kb \psi(x) |. 
$$

In the following, we assume that $f$ and $\sigma$ are smooth and 
since we are working with a stochastic differential equation on the torus, standard theorems give 
existence and uniqueness of a solution for any initial data $X(0)=x\in \T^d$. We denote this solution by 
$X(t,x),\; t\ge 0$. Also, since we chose to work on the torus, we do not have any problem of possible unbounded moments and this solution has 
clearly all moments finite.

We denote by $L(x,\partial_x)$ the Kolmogorov generator associated with the stochastic 
equation:
$$
L(x,\partial_x)v(x) = f^i(x) \partial_{i} v(x) + a^{ij}(x)\partial_{i j}v(x),
$$
where we use the summation convention for repeated indices and $\partial_i = \partial_{x_i}$ and 
$$
a^{ij}(x) := \frac12 \sum_{\ell = 1}^m \sigma^i_\ell(x) \sigma^j_\ell(x). 
$$

It is well known that the Kolmogorov equation:
\begin{equation}
\label{Eu0}
\frac{\dd u}{\dd t}=L(x,\partial_x)u,\; x\in \T^d,\; t>0,\quad u(0,x)=\varphi(x), \quad x\in \T^d,
\end{equation}
with periodic boundary conditions
has a unique solution for a smooth function $\varphi$ and that for all $x\in \T^d$:
$$
u(t,x)=\E(\varphi(X(t,x))).
$$
Moreover, this solution is smooth.
In the following, we write: $u(t)=P_t\varphi$ so that $(P_t)_{t\ge 0}$ is the transition semigroup 
associated with the Markov process $(X(t,x))_{t\ge 0,\; x\in \T^d}$. 
Note that we use the standard identification $u(t)=u(t,\cdot)$.

We wish to investigate the approximation properties of the Euler scheme for long times. We need
assumptions on the long time behavior of the law of the solutions of 
\eqref{Eu0}, {\it ie} of the law of the Markov process. We assume the following mixing properties:

\begin{itemize}
\item[\textbf{[H1]}] There exists a $\Ccal^\infty(\T^d,\R)$ function $\rho \geq 0$ such that 
\begin{equation}
\label{EH1}
L(x,\partial_x)^* \rho(x) = 0 \quad \mbox{and} \quad \int_{\T^d} \rho(x) \dd x = 1. 
\end{equation}
In other words, the measure $\rho(x) \dd x$ is invariant by $X(t,x)$. 
\item[\textbf{[H2]}]
Let $g \in \Ccal^\infty(\T^d,\R)$, and assume that $\int_{\T^d} g(x) \dd x = 0$. Then there exists a unique function $\mu(x) \in\Ccal^\infty(\T^d,\R)$ such that 
\begin{equation}
\label{EH2}
L(x,\partial_x)^* \mu(x) = g(x), \quad \mbox{and} \quad \int_{\T^d} \mu(x) \rho(x) \dd x = 0. 
 \end{equation}
\item[\textbf{[H3]}] Let $u(t,x)$ be the solution of \eqref{Eu0}. Assume that $\int_{\T^d} \varphi(x) \rho(x) 
\dd x  = 0$, then  there exists a constant $\lambda$ and, for each $k\in \N$, a polynomial $p_k(t)$,
 such that if $\varphi(x) \in \Ccal^\infty(\T^d,\R)$ we have the estimates 
\begin{equation}
\label{EH3}
\forall\, t \geq 0, \quad \Norm{u(t,\cdot)}{\Ccal^{k}} \leq p_k(t) e^{- \lambda t} \Norm{\varphi}{\Ccal^k}. 
\end{equation}
\end{itemize}

These hypothesis are usually satisfied under elliptic or hypoelliptic assumptions on the operator $L(x,\partial_x)$.  The reader
may find in  \cite{nualart} conditions to ensure {\bf [H1]}. We also refer to \cite{BallyTalay1, Mattingly2} for a general definition of hypoelliticity and applications to numerical schemes.  Combining kernel estimates for 
hypoelliptic diffusion (\cite{BallyTalay1, KusuokaStrook}) and exponential convergence to equilibrium, {\bf [H3]} can be proved. Note that similar estimates are used in \cite{Mattingly1, Talay90,Talay96,TalayLangevin}, where specific examples are considered.
Finally, we mention that  these hypothesis can be proved to be fulfilled using partial differential equations techniques (see \cite{hormander}).


For a smooth function $\psi$, we set 
$$
\langle \psi \rangle=\int_{\T^d} \psi(x)\rho(x)\dd x.
$$
Note that by {\bf [H3]}, we have for any solution of \eqref{Eu0} and $k\in \N$
\begin{equation}
\label{eq:exp0}
\forall\, t \geq 0, \quad \Norm{u(t,\cdot)-\langle \varphi\rangle}{\Ccal^{k}} \leq p_k(t) e^{- \lambda t} \Norm{\varphi - \langle \varphi\rangle}{\Ccal^k}. 
\end{equation}

Now for a small time step $\tau>0$ and $x\in \T^d$, we consider the Euler method defined, for $i = 1,\ldots, d$, by $X_0=x$ and the formula
\begin{equation}
\label{eq:Euler}
X^i_{n+1} = X^i_n + \tau f^i(X_n) + \sigma^i_\ell(X_n) (W^\ell((n+1)\tau)-W^\ell(n\tau)), 
\end{equation}
for $n \geq 0$. Our main result can be stated as follows: 
\begin{theorem}
Let $N$ and $\tau_0 > 0$ be fixed. Then there exists a modified smooth density 
$$
\mu^N(x) = \rho(x) + \tau \mu_1(x) + \cdots + \tau^N \mu_N(x) 
$$
such that $\int_{\T^d} \mu^N(x) \dd x =1$, a constant
$C_N$ and a polynomial $P_N(t)$ such that the following holds: For all smooth  function function $\varphi(x)$ on $\T^d$, we have 
\begin{equation}
\label{eq:approx3}
\forall\, p \in \N, \quad 
\Norm{\E\varphi(X_p) - \int_{\T^d} \varphi \dd \mu^N }{\infty} \leq \left(P_N(t_p) e^{- \lambda t_p } + C_N \tau^{N}
\right)\Norm{\varphi}{\Ccal^{8N + 2}}, 
\end{equation}
where for all $p$, $t_p = p \tau$ and $\dd \mu^N(x)  = \mu^N(x) \dd x$. 
\end{theorem}

This result can be viewed as a discrete version of \eqref{eq:exp0}. Note that it implies that all the invariant measure of the numerical process $X_p$ are close to $\dd \mu^N$ up to a very small error term $C_N \tau^N$. 

Using this result, we can also recover the weak convergence result
$$
\forall\, p \geq 1, \quad \Norm{\E\varphi(X_p) - \int_{\T^d} \varphi \dd \rho}{\infty} \leq C(\varphi) \tau
$$
for some constant $C(\varphi)$ depending on $\varphi$, and 
where we set $\dd \rho(x) = \rho(x) \dd x$. 
This can be compared with  \cite{Mattingly1, Mattingly2,Talay90, Talay96}. As in 
\cite{Mattingly2}, the only assumption made on $\tau$ is that $\tau \leq \tau_0$ where $\tau_0$ is any fixed number. The influence of $\tau_0$ is only reflected in the constants in the right-hand side - we can for example take $\tau_0 = 1$. In particular, we do not assume that $X_p$ has a unique invariant measure - something that would be guaranteed only if $\tau$ is small enough. We also recover an expansion of the invariant measure as in \cite{TalayTubaro}. 

In the next sections, the constants appearing in the estimate depend in general on bounds on derivatives of $f$ and $g$ defining the SDE. They will also depend in general on $\tau_0$ and $N$, but not on $\varphi$. 

\section{Asymptotic expansion of the weak error}

%
We 
have the formal expansion for small $t$: 
$$
u(t,x) = \varphi(x) + t L(x,\partial_x)\varphi(x) + \frac{t^2}{2} L(x,\partial_x)^2 \varphi(x) + \cdots + \frac{t^n}{n!} L(x,\partial_x)^n \varphi(x) + \cdots
$$ 
This is just obtained by Taylor expansion in time since by \eqref{Eu0} $\frac{\dd^n}{\dd t^n}u(t,x)=L(x,\partial_x)^nu(t,x)$ and in particular $\frac{\dd^n}{\dd t^n}u(0,x)=L(x,\partial_x)^n\varphi(x)$. 

Since the solution $u$ of the Kolmogorov equation is smooth and has its derivatives bounded 
in terms of the initial data $\varphi$, the above formal expansion can be 
justified and we have the following proposition whose proof is easy and left to the reader: 
\begin{proposition}
Assume that $\varphi \in \mathcal{C}^\infty(\T^d,\R)$, and let $\tau_0 > 0$. Then for all $N$, there exist a constant $C_N$ such that for all $\tau < \tau_0$, 
\begin{equation}
\label{EDA1}
\Norm{u(\tau,x) - \sum_{n = 0}^N  \frac{\tau^n}{n!}L(x,\partial_x)^n\varphi(x)  }{\infty} \leq C_N \tau^{N+1} \SNorm{\varphi}{\Ccal^{2N+2}}. 
\end{equation}
\end{proposition}

With the Euler scheme defined in \eqref{eq:Euler}, we associate the continuous process 
\begin{equation}
\label{EP1}
\tilde{X}_x^i(t) = X^i_n + (t-n\tau) f^i(X_n) + \sigma^i_\ell(X_n) (W^\ell(t)-W^\ell(n\tau)),\quad t\in [n\tau,(n+1)\tau],
\end{equation}
and $\tilde X(n\tau)=X_n$. We thus have $X_{n+1} = \tilde{X}_x(t_{n+1})$. 
The process \eqref{EP1} 
satisfies the equation
\begin{equation}
\label{EP2}
\dd \tilde{X}^i_x(t) = f^i(X_n)\dd t +   \sigma^i_\ell(X_n) \dd W^\ell(t),\quad t\in [n\tau,(n+1)\tau]. 
\end{equation}
Clearly, $(X_n)$ defines a discrete in time homogeneous Markov process but $\tilde X_x$ is not Markov.

In this work, we are only interested in the distributions of the solutions and of their approximation. 
We now examine in detail the first time step and its approximation properties in terms of the law. By Markov property, it is sufficient to then obtain information at all steps. Next result gives an expansion similar to 
Proposition \ref{EDA1} for the Euler process.
\begin{theorem} \label{t2.3}
Then for all $n\ge 1$, there exist operators $A_n(x,\partial_x)$ of order $2n$, such that for all $N\ge 1$, there exist a constant $C_N$ satisfying
\begin{equation}
\label{Eestd}
\Norm{\E\,\varphi(\tilde{X}_x(\tau)) - \sum_{n = 0}^N \tau^n A_n(x,\partial_x) \varphi(x)}{\infty} \leq C_N \tau^{N+1} \SNorm{\varphi}{\Ccal^{2N+2}},
\end{equation}
for all $\varphi \in \mathcal{C}^\infty(\T^d,\R)$ and $x\in \T^d$.
\end{theorem}
\begin{Proof}
Using \eqref{EP2} and the It\^o formula, we get for $t \leq \tau$,  
$$
\dd \varphi(\tilde X_x (\tau)) = L(x,\partial_x)\varphi(\tilde{X}_x(t)) + \sigma^i_\ell(x)  \partial_{i}\varphi(\tilde{X}_x(t))\dd W^\ell(x), 
$$
or equivalently,
\begin{equation}
\label{E1}
\varphi(\tilde X_x(t)) = \varphi(x) + \int_{0}^t L(x,\partial_x)\varphi(\tilde{X}_x(s)) \dd s + \int_{0}^t  \sigma^i_\ell(x) \partial_{i}\varphi(\tilde{X}_x(s)) \dd W^\ell(s). 
\end{equation}
Note that the last term is a martingale. We define the operator 
$$
R_{0,\ell}(x,\partial_x) = \sigma^i_\ell(x) \partial_{i}. 
$$
 We have for all $s$ and all $x$, 
$$
 L(x,\partial_x)\varphi(\tilde{X}_x(s)) = 
 (f^i(x) \partial_i \varphi)(\tilde{X}_x(s)) + (a^{ij}(x)\partial_{i j}\varphi)(\tilde{X}_x(s)). 
$$
Hence applying \eqref{E1} to $ \partial_i \varphi$ and $ \partial_{ij} \varphi$, we obtain
\begin{align*}
 L(x,\partial_x)\varphi(\tilde{X}_x(s))& = L(x,\partial_x) \varphi(x) \\[2ex]
& \displaystyle + \int_{0}^s  f^i(x) f^j(x) \partial_{ij} \varphi(\tilde X_x(\sigma))
 + f^i(x) a^{n\ell}(x)  \partial_{n\ell i} \varphi(\tilde X_x(\sigma)) \dd \sigma \\[2ex]
& \displaystyle  + \int_{0}^s 
  a^{ij}(x) f^n(x) \partial_{ijn} \varphi(\tilde X_x(\sigma))
 + a^{ij}(x) a^{n\ell}(x) \partial_{ijn\ell}  \varphi(\tilde X_x(\sigma))
 \dd \sigma\\[2ex]
& + \displaystyle \int_{0}^s  f^{i}(x) \sigma^j_\ell(x) \partial_{ij} \varphi(\tilde X_x(\sigma)) \dd W^\ell(\sigma)\\[2ex]
& + \displaystyle \int_{0}^s  a^{in}(x) \sigma^j_\ell(x) \partial_{inj} \varphi(\tilde X_x(\sigma)) \dd W^\ell(\sigma).
\end{align*}
We set $A_0=I$, $A_1=L$ and plugg this in \eqref{E1} to obtain
\begin{multline*}
\varphi(\tilde X_x(t)) =  \varphi(x) + t A_1(x,\partial_x) \varphi(x) + \int_{0}^t \int_{0}^s A_2(x,\partial_x) \varphi(\tilde X_x(\sigma)) \dd \sigma \dd s \\
  + \int_{0}^t \int_{0}^s R_{1,\ell}(x,\partial_x) \varphi(\tilde X_x(\sigma)) \dd W^\ell(\sigma) \dd s + \int_0^t R_{0,\ell}(x,\partial_x) \varphi (\tilde X_x(\sigma)) \dd W^\ell(\sigma). 
\end{multline*}
where 
$$
A_2(x,\partial_x) = f^i(x) f^j(x) \partial_{ij}  + f^i(x) a^{n\ell}(x)  \partial_{n\ell i}  + a^{ij}(x) f^n(x) \partial_{ijn} 
 + a^{ij}(x) a^{n\ell}(x) \partial_{ijn\ell}
$$
is an operator of order $4$, and 
$$
R_{1,\ell}(x,\partial_x) = f^{i}(x) \sigma^j_\ell(x) \partial_{ij}  + a^{in}(x) \sigma^j_\ell(x) \partial_{inj} 
$$
are operators of order $3$. Taking the expectation so that the last two term disappear, we
easily deduce the result for $N=1$.

Let us now prove recursively that 
there exist operators $A_n(x,\partial_x)$ of order $2n$ and $R_{n,\ell}(x,\partial_x)$ of order $2n+1$, such that 
\begin{equation}
\label{E2}
\begin{array}{rcl}
\varphi(\tilde{X}_x(t)) &=& \displaystyle \varphi(x) + t L(x,\partial_x) \varphi(x) + \sum_{n = 2}^{N} t^n A_n(x,\partial_x) \varphi(x) \\
&& \displaystyle + \int_{0}^t \cdots \int_{0}^{s_{N}} A_{N+1}(x,\partial_x)\varphi(\tilde{X}_x({s_{N+1}})) \dd s_1 \cdots \dd s_{N+1} \\[2ex]
&& + \displaystyle \sum_{n = 1}^{N} \int_{0}^t \cdots \int_{0}^{s_n} R_{n,\ell}(x,\partial_x)  \varphi(\tilde{X}_x({s_{n+1}})) \dd s_1 \cdots \dd s_{n} \dd W^\ell({s_{n+1}}). 
\end{array}
\end{equation}
Note the expectation of the last term vanishes so that \eqref{E2} easily implies \eqref{Eestd}. 

To prove \eqref{E2}, assume that  $A_{N+1}$ and $R_{N,\ell}$ are known, and let us decompose 
$A_{N+1}$ as  $A_{N+1}(x,\partial_x) = A^{\jb}_{N+1}(x) \partial_{\jb}$, where $\jb = (j_1,\ldots,j_m)$ are multiindices (with the summation convention) and $A^{\jb}_{N+1}$  smooth functions of $x$. 
Such decompistion is easy to write for $N=1$ or $2$. We apply \eqref{E1} to $\partial_{\jb}\varphi(\tilde{X}_x({s_{N+1}}))$ for a given muti-index $\jb$, and obtain
\begin{align*}
A^{\jb}_{N+1}(x) \partial_{\jb} \varphi(\tilde{X}_x({s_{N+1}})) &= A^{\jb}_{N+1}(x) \partial_{\jb} \varphi(x) \\[2ex]
&+ \displaystyle \int_{0}^{s_{N+1}} A^{\jb}_{N+1}(x) f^n(x) \partial_n \partial_{\jb} \varphi(\tilde{X}_x({s_{N+2}}))\dd s_{N+2} \\[2ex]
&+  \displaystyle\int_{0}^{s_{N+1}} A^{\jb}_{N+1}(x) a^{n\ell}(x) \partial_{n\ell} \partial_{\jb} \varphi(\tilde{X}_x({s_{N+2}}))\dd s_{N+2} \\[2ex]
& +\displaystyle \int_{0}^{s_{N+1}} A^{\jb}_{N+1}(x) \sigma_{\ell}^n(x) \partial_{n} \partial_{\jb} \varphi(\tilde{X}_x({s_{N+2}})) \dd W^\ell({s_{N+2}}). 
\end{align*}
We thus choose
$$
A_{N+2}(x) = A^{\jb}_{N+1}(x) f^n(x) \partial_n \partial_{\jb} + A^{\jb}_{N+1}(x) a^{ij}(x) \partial_{ij} \partial_{\jb}.
$$
and
$$
R_{N+1,\ell}(x) = A^{\jb}_{N+1}(x) \sigma_{\ell}^n(x) \partial_{n} \partial_{\jb},
$$
and obtain \eqref{E2} with $N+1$ replaced by $N+2$. 
\end{Proof}

\section{Modified generator}

\subsection{Formal series analysis}

Let  us now consider $\tau$ as fixed. We want to construct a formal series 
\begin{equation}
\label{e3.1}
L(\tau;x,\partial_x) =  L(x,\partial_x) + \tau L_1(x,\partial_x) + \cdots \tau^{n}L_n(x,\partial_x) + \cdots
\end{equation}
with operator coefficients $L_n(x,\partial_x)$ smooth on $\T^d$, and such  that formally the solution $v(t,x)$ at time $t = \tau$ of the equation
$$
\partial_t v(t,x) = L(\tau;x,\partial_x) v(t,x), \quad v(0,x) = \varphi(x)
$$
coincides in the sense of asymptotic expansion with the approximation of the transition 
semigroup $\E(\varphi(\tilde X_x(\tau)))$ studied in the previous section. In other words, we want to have the equality in the sense of asymptotic expansion in powers of $\tau$
$$
\exp( \tau L(\tau;x,\partial_x))\varphi(x) = \varphi(x) +  \sum_{n \geq 1} \tau^n A_n(x,\partial_x) \varphi(x), 
$$
where the operators $A_n(x,\partial_x)$ are defined in  Theorem \ref{t2.3}. 

Formally, this equation can be written
\begin{equation}
\label{Emlk}
\exp(\tau L(\tau;x,\partial_x)) - \mathrm{Id} = \tau \widetilde A(\tau) 
\end{equation}
where $\widetilde A(\tau) =  \sum_{n \geq 1} \tau^{n-1}  A_n$.

We have 
$$
\exp(\tau L(\tau;x,\partial_x)) - \mathrm{Id} = \tau L(\tau;x,\partial_x) \Big(\sum_{n \geq 0} \frac{\tau^n}{(n+1)!} L(\tau;x,\partial_x)^{n} \Big).
$$
Note that the (formal) inverse of the series is given by 
$$
\Big(\sum_{n \geq 0} \frac{\tau^n}{(n+1)!} L(\tau;x,\partial_x)^{n} \Big)^{-1} = \sum_{n \geq 0} \frac{B_n}{n!} \tau^n L(\tau;x,\partial_x)^n. 
$$
where the $B_n$ are the Bernoulli numbers: see for instance \cite{HLW,Fbook} and \cite{DF09} for a similar analysis involving operators. Hence equations \eqref{e3.1}, \eqref{Emlk} are equivalent in the sense of formal series to 
\begin{equation} 
\label{e3.3}
\begin{array}{lll}
L(\tau;x,\partial_x)\! &=&\! \ds \sum_{\ell \geq 0} \frac{B_\ell}{\ell!} \tau^\ell L(\tau;x,\partial_x)^\ell
\widetilde A(\tau)\\
\! &=& \! \ds \sum_{n \geq 0} \tau^n \left(A_{n+1} + \sum_{\ell = 1}^n \frac{B_\ell}{\ell!}  \sum_{n_1 + \cdots + n_\ell + n_{\ell+1} = n-\ell} L_{n_1} \cdots L_{n_\ell} A_{n_{\ell+1}+1}\right).
\end{array}
\end{equation}
Identifying the right hand sides of \eqref{e3.1} and \eqref{e3.3}, we get the following recursion formula
\begin{equation}
\label{eq:Ln}
L_n = A_{n +1}+ \sum_{\ell = 1}^n \frac{B_\ell}{\ell!}  \sum_{n_1 + \cdots + n_\ell + n_{\ell+1} = n-\ell} L_{n_1} \cdots L_{n_\ell} A_{n_{\ell+1}+1}. 
\end{equation}
Each of the terms of the above sum is an operator of order $2n+2$ with smooth coefficients and therefore  $L_n$ is also an operator of order $2n+2$ with smooth coefficients. 

Note that \eqref{Emlk} gives immediately the inverse relation of this formal series equation:
\begin{equation}
\label{eq:An}
A_n = \sum_{\ell = 1}^{n} \frac{1}{\ell!} \sum_{n_1 + \cdot + n_\ell = n - \ell} L_{n_1} \cdots L_{n_\ell} .
\end{equation}
Moreover, we have clearly
$$
L_n(x,\partial_x) \mathds{1} = 0,
$$
where $\mathds{1}$ denote the constant function equal to $1$.

\subsection{Approximate solution of the modified flow}

For a given $N$, we have constructed in the previous section an operator 
\begin{equation}
\label{ELN}
L^{(N)}(\tau;x,\partial_x) = L(x,\partial_x) + \sum_{n = 1}^N \tau^n L_n(x,\partial_x). 
\end{equation}
In order to perform  weak backward error analysis and estimate recursively the modified invariant law of the numerical process, we should be able to define a solution $v^N(t,x)$ of the modified flow 
\begin{equation}
\label{Emodif}
\partial_t v^N(t,x) = L^{(N)}(\tau; x,\partial_x) v^N(t,x), \quad v^N(0,x) = \varphi(x). 
\end{equation}
However, in our situation we do not know whether this equation has a solution. This is in contrast with standard backward error analysis where the modified flow can always be defined. 

The goal of the next proposition is to give a proper definition of the modified flow \eqref{Emodif}. 
%

\begin{theorem}
\label{t3.1}
Let $\varphi$ be a smooth functions on $\T^d$. For all $n \in \N$, there exist smooth functions $v_\ell(t,x)$, defined for all times $t \geq 0$, and such that 
for all $t \ge 0$ and $n \in \N$, 
\begin{equation}
\label{Ejkl}
\partial_t v_n(t,x) - L v_n(t,x) = \sum_{\ell = 1}^n {L_\ell} v_{n - \ell}(t,x),
\end{equation}
with initial conditions $v_0(0,x) = \varphi(x)$ and $v_n(0,x) = 0$ for $n \geq 1$. 
For all $N \geq 0$, setting 
\begin{equation}
\label{eq:vN}
\quad v^{(N)}(t,x) = \sum_{k = 0}^N \tau^N v_n(t,x), 
\end{equation} 
then the following holds: 

(i) There exists a constant $C_N$ such that for all time $t \geq 0$, and all $\tau \geq 0$, 
\begin{equation}
\label{eq:peniblos}
\Norm{\E v^{(N)}(t,\tilde X_x(\tau)) - v^{(N)}(t + \tau,x)}{\infty}\le C_N \tau^{N+1} \sup_{\substack{s \in (0,\tau) \\ n = 0, \ldots, N}}\SNorm{v_n(t+s)}{\Ccal^{4N + 2}}.  
\end{equation}

(ii) For fixed $\tau_0 > 0$, there exists a constant $C_N$ such that  for all $\tau \leq \tau_0$,  
$$
\Norm{\E\varphi(\tilde X_x(\tau)) - v^{(N)}(\tau,x)}{\infty}\le C_N \tau^{N+1} \Norm{\varphi}{\Ccal^{4N + 2}}. 
$$

\end{theorem}
\begin{Proof}For $n = 0$, the equation \eqref{Ejkl}  implies $v_0(t,x) = u(t,x)$, the solution of \eqref{Eu0}. 
Let $n\in\N$ and assume that $v_{j}(t,x)$ are constructed for $j=1\dots,n - 1$. Let 
\begin{equation}
\label{eq:Fn}
F_n(t,x) :=  \sum_{\ell = 1}^n {L_\ell} v_{n - \ell}(t,x)
\end{equation}
be the right-hand side in \eqref{Ejkl}.
Then $v_n$ is uniquely defined and given by the formula
\begin{equation}
\label{e4.3}
v_n(t)= \int_0^t P_{t-s}F_n(s)\dd s,\quad t\ge 0.
\end{equation}
By {\bf [H3]}, it is not 
difficult to check that $v_n,\; n\in \N$ are smooth and that for all $t \geq 0$, 
\begin{equation}
\label{eq:bapv}
\Norm{v_{n}(t)}{\Ccal^{k}}\le C(t) \Norm{\varphi}{\Ccal^{k+4n}},\quad k\in\N,
\end{equation}
where the constant $C(t)$ depends on $t$, $k$, $n$ and on the coefficient of the equation. Clearly, 
this type of estimate can be improved in the elliptic case. This proves the first part of the Theorem. 

To prove (i), we consider a fixed time $t$, and define the functions $w_n(x,s) := v_n(t + s)$ for $s \geq 0$ and $n \in \N$. 
By definition, these functions satisfy the relation
$$
\partial_s w_n(s,x) - L w_n(s,x) = \sum_{\ell = 1}^n {L_\ell} w_{n - \ell}(s,x), \quad w_n(0,x) = v_n(t,x). 
$$
Let us consider the successive time derivatives of the functions $w_n(s,x)$. We have using the definition of $w_n$, for all $s \geq 0$, 
$$
\partial_s^2 w_n(s,x) = \sum_{\ell = 0}^n L_\ell \partial_s w_{n - \ell}(s,x)
=\sum_{k = 0}^n \sum_{\ell_1 + \ell_2 = k} L_{\ell_1} L_{\ell_2} w_{n - k}(s,x), 
$$
and we see by induction that for all $m \geq 1$ and $s \geq 0$
\begin{equation}
\label{eq:ourf}
\partial_s^m w_n(s,x) =  \sum_{\ell_1 + \cdots + \ell_{m+1} = n} L_{\ell_1} \cdots L_{\ell_m} w_{\ell_{m+1}}(s,x). 
\end{equation}
Using  the fact that the operators $L_\ell$ are of order $\ell +2$ with no terms of order zero, we see that there exists a constant $C$ depending on $n$ and $m$, such that 
$$
\Norm{\partial_s^m w_n(s,x)}{\infty} \leq C \sup_{k = 0,\ldots,n}  \SNorm{w_{k}(s)}{\Ccal^{2k + 2m}}. 
$$

Now let us consider the Taylor expansion of $w_n(\tau)$, for $\tau \leq \tau_0$ and $n = 0,\ldots, N$, 
\begin{eqnarray*}
\label{eq:klm}
w_n(\tau,x) &=& \sum_{m = 0}^{N-n} \frac{\tau^m}{m! }\partial_t^m(0,x) + \int_{0}^\tau \frac{s^{N-m}}{(N-m)! }\partial_t^{N-m +1} w_n(s,x) \dd s  \\
&=& \sum_{m = 0}^{N-n} \frac{\tau^m}{m! }\sum_{\ell_1 + \cdots + \ell_{m+1} = n} L_{\ell_1} \cdots L_{\ell_m} w_{\ell_{m+1}}(0,x) + R_{N,n}(\tau,x). 
\end{eqnarray*}
Using the bounds on the time derivatives of $w_n(s,x)$, we obtain that for all $\tau \geq 0$ and  all $n = 0, \ldots, N$, 
$$
\Norm{R_{N,n}(\tau,x)}{\infty}
\leq C \tau^{N - m + 1} \sup_{\substack{s \in (0,\tau) \\ n = 0, \ldots, N}} \SNorm{w_n(s,x)}{\Ccal^{4N + 2}}, 
$$
for some constant depending on $N$, $m$. After summation in $n$, and using the expression \eqref{eq:An} of the operators $A_n$ and the definition of $w_n$, we get
$$
v^{(N)}(t + \tau,x) = \sum_{n = 0}^N \tau^n \sum_{m = 0}^{n} A_m v_{n - m}(t,x)  
+ R_N(t,\tau, x)
$$
where 
$$
\Norm{R_N(t,\tau, x)}{\infty} \leq C_N \tau^{N + 1} \sup_{\substack{s \in (0,\tau) \\ n = 0, \ldots, N}} \SNorm{v_n(t+s,x)}{\Ccal^{4N + 2}}.
$$
To conclude, we use \eqref{Eestd} applied to $\varphi = v^{(N)}(t,x)$, and we easily verify that $\E v^{(N)}(\tau, \tilde X_x(\tau))$ satisfy the same asymptotic expansion. 

The second estimate (ii) is then a consequence of (i) with $t = 0$ and \eqref{eq:bapv}. 
\end{Proof}

Note that in the previous theorem, we have constructed a function $v^{(N)}(t,x)$ which is an approximate solution of \eqref{Emodif}. More precisely, we can easily show that we have for all time $t \geq 0$, 
$$
\partial_t v^{(N)}(t,x)=L^{(N)}(\tau;x,\partial_x)v^{(N)}(t,x) +R^{(N)}(t,x),\quad v^{(N)}(0,x)=\varphi(x),
$$
where 
$$
R^{(N)}(t,x) = - \sum_{\substack{\ell_1,\ell_2=0,\dots,N \\ \ell_1+\ell_2>N }} \tau^{\ell_1+\ell_2} L_{\ell_1}
v_{\ell_2}(t,x)
$$
is of order $\mathcal{O}(\tau^{N+1})$.

\section{Asymptotic expansion of the invariant measure and long time behavior}

We now analyze the long time behavior of the solution of the modified equation \eqref{Emodif}. In the following, for a given operator $B(x,\partial_x)$, we denote by  $B(x,\partial_x)^*$ its adjoint with respect to the $L^2$ product. We start by an asymptotic expansion of a formal invariant measure
for the numerical scheme.

\begin{proposition}
Let $(L_n)_{n \geq 0}$ be the collection of operators defined recursively by \eqref{eq:Ln}. There exists a collection of functions $(\mu_n(x))_{n \geq 0}$ such that $\mu_0(x) = \rho(x)$, $\int_{\T^d} \mu_n(x) \dd x = 0$ for $n \geq 1$, and for all $n \geq 0$, 
\begin{equation}
\label{Emuk}
L_0^{*} \mu_n = - \sum_{\ell =1}^{n} (L_\ell)^* \mu_{n-\ell}. 
\end{equation}
Let $N \geq 0$ be fixed and $L^{(N)}(\tau;x,\partial_x)$ the operator defined by \eqref{ELN}. Then the function 
$$
\mu^{(N)}(\tau;x) = \rho(x) + \sum_{n = 1}^N \tau^n \mu_n(x) \in \Ccal^\infty(\T^d,\R)
$$
satisfies 
$$
\int_{\T^d} \mu^{(N)}(\tau;x) \dd x = 1, 
$$
and 
$$
L^{(N)}(\tau;x,\partial_x)^* \mu^{(N)}(\tau;x) = G^{(N)}(\tau;x),  
$$
with for all $k$ and all $\tau \in [0,\tau_0]$, 
$$
\Norm{G^{(N)}(\tau;x)}{\Ccal^k} \leq C_{N,k} \tau^{N+1} \quad \mbox{and}\quad \int_{\T^d} G^{(N)}(\tau;x) \dd x = 0, 
$$
where $C_{N,k}$ depends on $N$ and $k$. 
\end{proposition}
\begin{Proof}
Assume that $\mu_0 = \rho$ and $\mu_j(x)$ are known, for $j = 0,\ldots, n-1$ with $n \geq 1$. 
Consider the equation \eqref{Emuk} given by 
$$
L_0^{*} \mu_n = - \sum_{\ell =1}^{n} (L_\ell)^* \mu_{n-\ell} = G_n. 
$$
Note that the right-hand side $G_n(x)$ is a smooth function satisfying 
$$
\int_{\T^d} G_n(x) \dd x = - \sum_{\ell =1}^{n} \int_{\T^d}(L_\ell)^* \mu_{n-\ell} \dd x =  - \sum_{\ell =1}^{n} \int_{\T^d}\mu_{n-\ell}L_\ell \mathds{1} \dd x = 0, 
$$
where $\mathds{1}$ denotes the constant function equal to $1$, which as already seen is in the Kernel of all the $L_\ell$.

Using Hypothesis \textbf{[H2]}, we easily obtain the existence of a $\Ccal^\infty$ function $\mu_n$ satisfying \eqref{Emuk} and $\int_{\T^d} \mu_n(x) \dd x = 0$. This shows the first part of the Proposition. 

We then write  
\begin{eqnarray*}
L^{(N)}(\tau;x,\partial_x)^* \mu^{(N)}(\tau;x) &=& \sum_{n = 0}^{2N} \tau^n\sum_{\substack{\ell_1 + \ell_2 = n \\\ell_i \leq N}} (L_{\ell_1})^* \mu_{\ell_2}\\
&=& \sum_{n = N+1}^{2N} \tau^n\sum_{\substack{\ell_1 + \ell_2 = n \\ \ell_i \leq N}} (L_{\ell_1})^* \mu_{\ell_2}\\
&=:& G^{(N)}(\tau;x), 
\end{eqnarray*}
and we easily verify that $G^{(N)}$ satisfies the hypothesis of the Proposition, owing to the fact that $L_\ell\mathds{1} = 0$ for all $\ell \geq 0$. 
\end{Proof}

\begin{proposition}
\label{p4.2}
For all $n$ and $k$ there exists a polynomial $P_{k,n}(t)$ such that for all $t \geq 0$, 
\begin{equation}
\label{eq:exp}
\Norm{v_n(t,x) - \int_{\T^d} \varphi  \dd \mu_n }{\Ccal^k} \leq  P_{k,n}(t)  e^{-\lambda t} \Norm{\varphi - \langle \varphi \rangle}{\Ccal^{k + 4n}}.
\end{equation}
\end{proposition}

\begin{Proof}
Using the fact that $\mu_0 = \rho$ and $v_0 = u$, wee see that estimate \eqref{eq:exp} is satisfied for $n = 0$ (see Equation \eqref{eq:exp0}). 
Let $n\geq 1$ and assume that $v_j,\; j=0,\dots, n-1$ satisfy for $k\in\N$, $t\ge 0$:
$$
\Norm{v_j(t,x) - \int_{\T^d}\varphi(x) \mu_j(x) \dd x}{\Ccal^k} \leq P_{k,j}(t) e^{-t\lambda}\Norm{\varphi(x) - \langle \varphi\rangle}{\Ccal^k},
$$
for some polynomial $P_{k,j}$.

Let us set 
$$
c_n(t) = \sum_{m = 0}^n \int_{\T^d} v_{n-m}(t,x)  \mu_m(x) \dd x.
$$
We claim that $c_n(t)$ does not depend on time. Indeed:
\begin{eqnarray*}
\sum_{m = 0}^n \partial_t \int_{\T^d} v_{n-m}(t,x)  \mu_m(x) \dd x  &=& \sum_{m = 0}^n \partial_t \int_{\T^d} v_{m}(t,x)  \mu_{n-m}(x) \dd x\\
&=&\sum_{m=0}^n\sum_{\ell=0}^m  \int_{\T^d} L_{m-\ell}v_{\ell}(t,x) \mu_{n-m}(x) \dd x\\
&=& \sum_{\ell=0}^n \sum_{m=\ell}^n  \int_{\T^d}  v_{\ell}(t,x)  L_{m-\ell}^* \mu_{n-m}(x) \dd x\\
&=& \sum_{\ell=0}^n \int_{\T^d} v_{\ell}(t,x)   \sum_{m=0}^{n-\ell} L_{m}^* \mu_{n-\ell-m}(x) \dd x  = 0, 
\end{eqnarray*}
by definition of the coefficients $\mu_n$, see \eqref{Emuk}. Note that, thanks to the smoothness
properties of all the functions, the computation above is easily justified. 

We deduce:
\begin{equation}
\label{e4.4}
 \int_{\T^d} \partial_t v_n(t,x) \rho(x) \dd x = - \sum_{m = 1}^n \int_{\T^d} \partial_t v_{n-m}(t,x) \mu_m(x) \dd x.
\end{equation}

Next, we compute the average of $F_n$. By \eqref{Ejkl}, \eqref{eq:Fn} and  \eqref{e4.4}, we have 
\begin{eqnarray*}
\langle F_n(t)\rangle=\int_{\T^d} F_n(t,x) \rho(x) \dd x &=& \int_{\T^d} \partial_t v_n(t,x) \rho(x) \dd x -\int_{\T^d} Lv_n(t,x)\rho(x)\dd x\\
&=&  \int_{\T^d} \partial_t v_n(t,x) \rho(x) \dd x\\
&=& -\sum_{m = 1}^n \int_{\T^d} \partial_t v_{n-m}(t,x) \mu_m(x) \dd x. 
\end{eqnarray*}

We rewrite \eqref{e4.3} as follows
$$
v_n(t,x) = \int_0^t \langle F_n(s) \rangle \dd s + \int_{0}^t P_{t-s} (F_n(s,x) - \langle F_n(s) \rangle) \dd s.
$$ 
Using the previous expression obtained for $\langle F_n(s) \rangle$ and recalling the initial data 
for $v_n$, we deduce that 
\begin{multline*}
v_n(t,x) =  
-\sum_{m = 1}^n \int_{\T^d} v_{n-m}(t,x) \mu_m(x) \dd x + \int_{\T^d} \varphi(x) \mu_n(x) \dd x \\ + 
\int_{0}^t P_{t-s} (F_n(s,x) - \langle F_n(s) \rangle) \dd s. 
\end{multline*}
Then, using  $\int_{\T^d} \mu_m(x) \dd x = 0$, $m\in \N$, we get 
\begin{multline*}
v_n(t,x) - \int_{\T^d} \varphi(x) \mu_n(x) \dd x = \sum_{m = 1}^n \int_{\T^d} (v_{n-m}(t,x) - \int_{\T^d} \varphi(x) \mu_{n-m}(x) \dd x) \mu_m(x) \dd x \\
+ \int_{0}^t P_{t-s} (F_n(s,x) - \langle F_n(s) \rangle) \dd s. 
\end{multline*}
Note that, since $L_\ell$, $\ell\in \N$ is a differential operator of order $2\ell+2$ with smooth 
coefficients and containing no zero order terms, we have
$$
\Norm{F_n(s) - \langle F_n(s) \rangle}{\Ccal^k} \le c_{k,\ell} \sum_{\ell=0}^{n-1}\Norm{v_{\ell}(s)
-\int_{\T^d}v_\ell(s,x) \mu_\ell(x)dx}{\Ccal^{k+2(n-\ell)+2}}. 
$$
Then by {\bf [H3]}
\begin{multline*}
\Norm{ v_n(t,x) - \int_{\T^d} \varphi(x) \mu_n(x) \dd x}{\Ccal^k} \le \\ \sum_{m = 1}^n \Norm{v_{n-m}(t,x) - \int_{\T^d} \varphi(x) \mu_{n-m}(x) \dd x}{\infty}  \int_{\T^d}|\mu_m(x)| \dd x \\
+ \int_{0}^t p_k(t-s) e^{-\lambda (t-s)} \Norm{F_n(s,x) - \langle F_n  (s) \rangle}{\Ccal^k} \dd s,\\
\end{multline*}
and using the recursion assumption
\begin{multline*}
\Norm{ v_n(t,x) - \int_{\T^d} \varphi(x) \mu_n(x) \dd x}{\Ccal^k} \ds \le  \sum_{m = 1}^n c_{n,m}P_{0,n-m}(t) e^{-t\lambda}\Norm{\varphi(x)-\langle \varphi\rangle}{\infty}\\
\ds + \sum_{\ell=0}^{n-1}\int_{0}^t p_k(t-s) P_{k+2(n-\ell)+2, \ell}(s) e^{-\lambda t} \dd s \Norm{\varphi(x)-\langle \varphi\rangle}{k+4n}. 
\end{multline*}
The conclusion follows.
\end{Proof}

%

We give now our main result concerning the long time behavior of the numerical solution: 

\begin{theorem}
Let $\tau_0$ and $N$ be fixed. Then there exists $C_N$ and a polynomial $P_N(t)$ such that the following holds: 
Let $X_p$ be the discrete process defined by \eqref{eq:Euler}, then we have for $p \geq 0$, $\tau\geq 0$ and smooth function $\varphi(x)$
\begin{equation}
\label{eq:approx}
\forall\, p \in \N, \quad 
\Norm{\E \, \varphi(X_p) - v^{(N)}(t_p ,x)}{\infty} \leq C_N\tau^{N} \Norm{\varphi}{\Ccal^{8N+2}},
\end{equation}
where for all $p$, $t_p = p \tau$. Moreover, we have 
\begin{equation}
\label{eq:approx2}
\forall\, p \in \N, \quad 
\Norm{\E\varphi(X_p) - \int_{\T^d} \varphi \dd \mu^{(N)} }{\infty} \leq \left(P_N(t_p) e^{- \lambda t_p } + C_N \tau^{N}, 
\right)\Norm{\varphi}{\Ccal^{8N + 2}}, 
\end{equation}
where $\dd \mu^{(N)}(x) = \mu^{(N)}(x) \dd x$. 
\end{theorem}
\begin{Proof}
For all $p$, with $t_j = j\tau$, we have 
\begin{eqnarray*}
\E \, \varphi(X_p) - v^{(N)}(t_p,x) &=& \E\, v^{(N)}(0,X_p) - v^{(N)}(t_p,x)\\
&=& \E \sum_{j = 0}^{p-1} \E^{X_{p-j-1}} \Big(v^{(N)}(t_j,X_{p - j}) - v^{(N)}(t_{j+1}, X_{p - j - 1}) \Big). 
\end{eqnarray*}
Here we have used  the notation 
$ \E^{X_{p-j-1}} $ for the conditional expectation with respect to the filtration generated by $X_{p-j-1}$.
By the Markov property of the Euler process at times $t_j$:
\begin{multline*}
\ds \E^{X_{p-j-1}} \Big(v^{(N)}(t_j,X_{p - j}) - v^{(N)}(t_{j+1}, X_{p - j - 1}) \Big)\\
\ds =\E^{X_{p-j-1}} \Big(v^{(N)}(t_j,\tilde X_{X_{p - j - 1}}(\tau)) - v^{(N)}(t_{j+1}, X_{p - j - 1}) \Big).
\end{multline*}
Using \eqref{eq:peniblos} with $t = t_j$, 
and Proposition \ref{p4.2}, we deduce that 
\begin{eqnarray*}
\Norm{\E \, \varphi(X_p) - v^{(N)}(t_p,x)}{\infty} &\le & C_N \tau^{N+1}\sum_{j=0}^{p-1} \sup_{\substack{s \in (0,\tau) \\ n = 0,\ldots, N}} \SNorm{v_n(t_j + s, x)}{4N + 2} \\ 
 &\le &  C_N\tau^{N+1}\Norm{\varphi}{\Ccal^{8N + 2}}  \sum_{j=0}^{p-1}Q_N(t_j)
e^{-\lambda t_j} 
\end{eqnarray*}
for some constant $C_N$ and polynomial $Q_N(t)$.  We have used:
$\SNorm{v_n(t_j + s, x)}{4N + 2}=\SNorm{v_n(t_j + s, x)-\int_{\T^d}\varphi d\mu_n}{4N + 2}$.
We conclude by using the fact that for a fixed constant $\gamma > 0$, we have 
$$
\sum_{j = 0}^{p-1} e^{- \gamma j \tau} \leq \frac{1}{ 1 - e^{- \gamma \tau}} \leq \frac{C}{\tau}
$$
where the constant $C$ depends on $\gamma$ and $\tau_0$. 
This shows \eqref{eq:approx}. The second estimate is a consequence of 
Proposition \ref{p4.2}. 
\end{Proof}

\end{document}